\documentclass[12pt]{amsart}
\usepackage[centertags]{amsmath}
\usepackage{amsfonts,amssymb}
\usepackage{graphicx}
\usepackage{enumerate}
\usepackage[latin1]{inputenc}
\usepackage{float}

  \newtheorem{theorem}{Theorem}
  \newtheorem{corollary}{Corollary}

  \theoremstyle{remark}
  \newtheorem{remark}{Remark}

 \newcommand{\lcm}{{\rm lcm}}
 \newcommand{\dd}{{\rm d}}
 \newcommand{\ddiv}{{\rm div}}

\begin{document}

\title{Permutations with arithmetic constraints}
\date{\today}

\author{Carl Pomerance}
\address{Mathematics Department, Dartmouth College, Hanover, NH 03784}
\email{carlp@math.dartmouth.edu}

\begin{abstract}
Let $S_\lcm(n)$ denote the set of permutations $\pi$ of $[n]=\{1,2,\dots,n\}$
such that $\lcm[j,\pi(j)]\le n$ for each $j\in[n]$.  Further, let $S_\ddiv(n)$ denote
the number of permutations $\pi$ of $[n]$ such that
$j\mid\pi(j)$ or $\pi(j)\mid j$ for each $j\in[n]$.  Clearly $S_\ddiv(n)\subset S_\lcm(n)$.
We get upper and lower bounds for the counts of these sets, showing they 
grow geometrically.  We also prove a conjecture from a recent paper on the number
of ``anti-coprime" permutations of $[n]$, meaning that each $\gcd(j,\pi(j))>1$ except
when $j=1$.
\end{abstract}


\maketitle

\vskip-30pt
\newenvironment{dedication}
       {\vspace{6ex}\begin{quotation}\begin{center}\begin{em}}
       {\par\end{em}\end{center}\end{quotation}}
\begin{dedication}
{In memory of Eduard Wirsing (1931--2022)}
\end{dedication}
\vskip20pt

\section{Introduction}

Recently, in \cite{P22} some permutation enumeration problems with an arithmetic
flavor were considered.  In particular, one might count permutations $\pi$ of $[n]=\{1,2,\dots,n\}$
where each $\gcd(j,\pi(j))=1$ and also permutations $\pi$ where each $\gcd(j,\pi(j))>1$
except for $j=1$.  It was shown in \cite{P22} that the coprime count is between
$n!/c_1^n$ and $n!/c_2^n$ for all large $n$, where $c_1=3.73$ and $c_2=2.5$.
Shortly after, Sah and Sawhney \cite{SS} showed that there is an explicit constant $c_0=2.65044\dots$
with the count of the shape $n!/(c_0+o(1))^n$ as $n\to\infty$.  The ``anti-coprime" count
was shown in \cite{P22} to exceed $n!/(\log n)^{(\alpha+o(1))n}$ as $n\to\infty$, where
$\alpha=e^{-\gamma}$, with $\gamma$ Euler's constant.  
It was conjectured in \cite{P22} that this lower bound is sharp, which we will prove here.

There are several papers in the literature that have considered the divisibility graph
on $[n]$ where $i\ne j$ are connected by an edge if $i$ divides $j$ or vice versa,
and the closely related lcm graph, where edges correspond to $\lcm[i,j]\le n$.
In particular, it was shown in \cite{P2} that the length of the longest simple
path in such graphs is $o(n)$, and this has been improved to order-of-magnitude $n/\log n$,
see Saias \cite{S} for a recent paper on the topic.  One might also consider permutations
of $[n]$ compatible with these graphs.
Let  $S_\ddiv(n)$ denote the set of permutations
$\pi$ of $n$ such that for each $j\in[n]$, either $j\mid\pi(j)$ or $\pi(j)\mid j$.
Further, let
$S_{\lcm}(n)$ denote the set of permutations $\pi$ of $[n]$ such
that for each $j\in[n]$, $\lcm[j,\pi(j)]\le n$.   Clearly, $S_\ddiv(n)\subset S_\lcm(n)$.
There is a small literature on these topics.  In particular, counts for $\#S_\ddiv(n)$
are on OEIS \cite{oeis} (due to Heinz and Farrokhi), which we reproduce here, together with new counts for $\#S_\lcm(n)$.

\begin{footnotesize}
\begin{table}[]
\caption{Counts for $S_\ddiv(n)$ and $S_\lcm(n)$ and their $n$th roots.}
\label{Ta:counts}
\begin{tabular}{|rrcrcc|} \hline
$n$&$\#S_\ddiv(n)$ & $(\#S_\ddiv(n))^{1/n}$&$\#S_\lcm(n)$&$(\#S_\lcm(n))^{1/n}$&\\ \hline
1	&	1&1.0000&1&1.0000&\\
2	&	2&1.4142&2&1.4142&\\
3	&	3&1.4422&3&1.4422&\\
4	&	8&1.6818&8&1.6818&\\
5	&	10&1.5849&10&1.5849&\\
6	&	36&2.8272&56&1.9560&\\
7	&	41&1.6998&64&1.8114&\\
8	&	132&1.8411&192&1.9294&\\
9	&	250&1.8469&332&1.9060&\\
10	&	700&1.9254&1{,}184&2.0292&\\
11	&	750&1.8254&1{,}264&1.9142&\\
12	&	4{,}010&1.9965&12{,}192&2.1903&\\
13	&	4{,}237&1.9011&12{,}872&2.0708&\\
14	&	10{,}680&1.9398&37{,}568&2.1221&\\
15	&	24{,}679&1.9626&100{,}836&2.1556&\\
16	&	87{,}328&2.0362&311{,}760&2.2048&\\
17	&	90{,}478&1.9569&322{,}320&2.1087&\\
18	&	435{,}812&2.0573&2{,}338{,}368&2.2585&\\
19	&	449{,}586&1.9839&2{,}408{,}848&2.1671&\\
20	&	1{,}939{,}684&2.0625&14{,}433{,}408&2.2802&\\
21	&	3{,}853{,}278&2.0588&32{,}058{,}912&2.2773&\\
22	&	8{,}650{,}900&2.0669&76{,}931{,}008&2.2828&\\
23	&	8{,}840{,}110&2.0046&78{,}528{,}704&2.2043&\\
24	&	60{,}035{,}322&2.1091&919{,}469{,}408&2.3631&\\
25	&	80{,}605{,}209&2.0714&1{,}158{,}792{,}224&2.3044&\\
26	&	177{,}211{,}024&2.0761&2{,}689{,}828{,}672&2.3051&\\
27	&	368{,}759{,}752&2.0757&4{,}675{,}217{,}824&2.2811&\\
28	&	1{,}380{,}348{,}224&2.1205&21{,}679{,}173{,}184&2.3396&\\
29	&	1{,}401{,}414{,}640&2.0673&21{,}984{,}820{,}864&2.2731&\\
30	&	8{,}892{,}787{,}136&2.1460&381{,}078{,}324{,}992&2.4324&\\
31	&	9{,}014{,}369{,}784&2.0947&386{,}159{,}441{,}600&2.3646&\\
32	&	33{,}923{,}638{,}848&2.1334&1{,}202{,}247{,}415{,}040&2.3851&\\
33	&	59{,}455{,}553{,}072&2.1208&&&\\
34	&	126{,}536{,}289{,}568&2.1210&&&\\
35	&	207{,}587{,}882{,}368&2.1055&&&\\
\hline
\end{tabular}
\end{table}
\end{footnotesize}

The table suggests that $\#S_\lcm(n)>\#S_\ddiv(n)>2^n$ for $n$ large,
and that there may be a similar upper bound.  
In this note we will prove 
that $(\#S_\ddiv(n))^{1/n}$ is bounded above 1 and $(\#S_\lcm(n))^{1/n}$
is bounded below infinity.  We conjecture they tend to limits, but we lack the
numerical evidence or heuristics to suggest  values for these limits.
 We will also show that $\#S_\lcm(n)/\#S_\ddiv(n)$ tends to
infinity geometrically.

One might also ask for the length of the longest cycle among permutations in
$S_\ddiv(n)$ or in $S_\lcm(n)$.  This seems to be only slightly less (if at all)
than the length of the longest simple chain in the divisor graph or lcm graph on $[n]$ 
mentioned above.
Other papers have looked at tilings of $[n]$ with divisor chains, for example see
\cite{M}.  This could correspond to asking about the cycle decomposition for
permutations in $S_\ddiv(n)$ or in $S_\lcm(n)$.

We mention the paper \cite{EFH} of Erd\H os, Freud, and Hegyv\'ari where
some other arithmetic problems connected with integer permutations are discussed.
Finally, we note the recent paper \cite{BB} which also has a similar flavor.

\section{An upper bound for $\#S_{\lcm}(n)$}

\begin{theorem}
\label{main}
We have $\# S_{\lcm}(n)\le e^{2.61n}$ for all large $n$.
\end{theorem}
\begin{proof}
Let $n$ be large.  For $j\in[n]$, let $N(j)$ denote the number of $j'\in[n]$ with
$\lcm[j,j']\le n$.  This condition can be broken down as follows:  $\lcm[j,j']\le n$
if and only if there are integers $a,b,c$ with
\begin{equation}
\label{eq:abc}
 j=ab,\quad j'=bc,\quad \gcd(a,c)=1,\quad abc\le n.
\end{equation}
That is, $N(j)$ is the number of triples $a,b,c$ with $j=ab$ satisfying \eqref{eq:abc}.
Since $a\mid j$, $b=j/a$, and $c\le n/j$, we have $N(j)\le\tau(j)n/j$, where $\tau$ is the divisor function 
(which counts the number of positive divisors of its argument). 
For $\pi\in S_\lcm(n)$, the number of possible values for $\pi(j)$ is at most $N(j)$, so we have
\begin{equation}
\label{eq:Nj}
\#S_\lcm(n)\le\prod_{j\in[n]}N(j)\le \prod_{j\in[n]}\tau(j)n/j.
\end{equation}
This quickly leads to an estimate for
$\#S_\lcm(n)$ that is of the form $n!^{o(1)}$ as $n\to\infty$, but to do better we
will need to work harder.  In particular we use a seemingly trivial property of
permutations: they are one-to-one.  In particular, there are not many values of
$j$ with $\pi(j)$ small, since there are not many small numbers.  This thought leads
to versions of \eqref{eq:Nj} where $\tau$ is replaced with a restricted divisor function
that counts only small divisors.

Let $k=30$.  We partition the interval $(0,n]$ into
subintervals as follows.  Let $J_0=(n/k,n]$.  Let $i_0$ be the largest $i$ such that
$L:=k^{2^i}\le\log n$, so that $(\log n)^{1/2}<L\le \log n$.  For $i=1,\dots,i_0$, let
$J_i=(n/k^{2^i},n/k^{2^{i-1}}]$, and let $J_{i_0+1}=(0,n/L]$.

For $\pi\in S_\lcm(n)$ we have sets $X_i,Y_i$ as follows:
\begin{align*}
X_i:=&\{j\in J_i:\pi(j)> n/k^{2^i}\},\quad 0\le i\le i_{0},\\
Y_i:=&\{j>n/k^{2^{i-1}}:\pi(j)\in J_i\},~1\le i\le i_{0}+1.
\end{align*}
These sets depend on the choice of $\pi$, but the number of
choices for the sets $Y_i$ is not so large.  We begin by counting the number of
possibilities for the sequence of sets $Y_1,\dots,Y_{i_0+1}$.

Since $\pi$ is a permutation it follows that $y_i:=\#Y_i$ is at most the number
of integers in $J_i$, so that $y_i\le n/k^{2^{i-1}}$.  The number of subsets of
$(n/k^{2^{i-1}},n]$ of cardinality $\le y_i$ is less than
\[
\sum_{u\le y_i}\binom nu \le 2\binom n{y_i}
\le\frac{2n^{y_i}}{y_i!}
\le\exp\Big(\frac n{k^{2^{i-1}}}(2^{i-1}\log k+1)\Big),
\]
for $n$ sufficiently large, using the inequality $2/j!<(e/j)^j$ for $j\ge3$.  Multiplying these estimates we obtain that the number
of choices for a sequence of sets $\{Y_i\}$ as described is
\begin{equation}
\label{eq:Yseq}
\le\exp(0.1554n)
\end{equation}
for all sufficiently large $n$.

Fix now a specific sequence of sets $\{Y_i\}$, which then determines a complementary
sequence of sets $\{X_i\}$ with $X_i=J_i\setminus Y_{i+1}$.  We will give the set $X_0$ special treatment, so for now,
assume that $1\le i\le i_0$.  For $j\in X_i$, the number of possible choices $j'$ to which $j$ may be
mapped by a permutation in $S_\lcm(n)$ (with sequence of sets $\{Y_i\}$) is at most
the number of choices for $a,c$ as in \eqref{eq:abc}.  Here $c\le n/j$ and $a\mid j$
with $a\le n/j'<k^{2^i}$.  Let $\tau_z(m)$ be the number of divisors of $m$ that are
$<z$.  With this notation, the number of choices for $\pi\in S_\lcm(n)$ restricted to $X_i$ is at most
\begin{align*}
\prod_{j\in X_i}\tau_{k^{2^i}}(j)\frac nj&\le\prod_{j\le n/k^{2^{i-1}}}\tau_{k^{2^i}}(j)\frac nj\\
&\le\left(\frac1{\lfloor n/k^{2^{i-1}}\rfloor}\sum_{j\le n/k^{2^{i-1}}}\tau_{k^{2^i}}(j)\right)^{n/k^{2^{i-1}}}\frac{n^{n/k^{2^{i-1}}}}{\lfloor n/k^{2^{i-1}}\rfloor!},
\end{align*}
by the AM-GM inequality (the arithmetic mean geometric mean inequality).

Since the harmonic sum $\sum_{d<z}1/d$ is bounded above by $\log z+1$, we have
\begin{equation}
\label{eq:tauz}
\sum_{j\le x}\tau_z(j)=\sum_{d<z}\left\lfloor\frac xd\right\rfloor
=\sum_{d<z}\left\lfloor\frac{\lfloor x\rfloor}d\right\rfloor<\lfloor x\rfloor(\log z+1).
\end{equation}
We apply this above getting that the number of assignments for the numbers $j\in X_i$ is
at most
\begin{align*}
\big(\log(k^{2^i})+1\big)^{n/k^{2^{i-1}}}&\exp\Big(\frac n{k^{2^{i-1}}}(\log(k^{2^{i-1}})+1)\Big)\\
=&\exp\Big(\frac n{k^{2^{i-1}}}\big(\log(k^{2^{i-1}})+\log(\log(k^{2^i})+1)+1\big)\Big).
\end{align*}
Thus, multiplying these estimates we have that the number of assignments for numbers $j$
in the sets $X_i$, $1\le i\le i_0$ is
\begin{equation}
\label{eq:Xest}
\le\exp(0.2269n)
\end{equation}
for $n$ sufficiently large.

We next deal with the elements of the sets $Y_i$.  Again referring to \eqref{eq:abc}, for each
$j\in Y_i$ we are to count pairs $a,c$ with $a\mid j$ and $a\le n/j'<k^{2^i}$ and $c\le n/j<k^{2^{i-1}}$.
Assuming $Y_i$ is not empty, the number of assignments for elements of $Y_i$ is at most
\begin{align*}
\prod_{j\in Y_i}\big(\tau_{k^{2^i}}(j)k^{2^{i-1}}\big)&
\le\left(\frac1{y_i}\sum_{j\in Y_i}\tau_{k^{2^i}}(j)k^{2^{i-1}}\right)^{y_i}\\
&\le\left(\frac1{y_i}\sum_{j\le n}\tau_{k^{2^i}}(j)k^{2^{i-1}}\right)^{y_i}\\
&\le\left(\frac n{y_i}(\log(k^{2^i})+1)k^{2^{i-1}}\right)^{y_i},
\end{align*}
using \eqref{eq:tauz}.  We have $y_i<n/k^{2^{i-1}}$ and in this range, the above estimate
is increasing as $y_i$ varies.  So, the count is at most
\[
\exp\Big(\frac n{k^{2^{i-1}}}\big(\log(k^{2^i})+\log(\log(k^{2^i})+1)\big)\Big).
\]
Multiplying these estimates we have that the number of assignments for numbers $j$ in
the sets $Y_i$ is
\begin{equation}
\label{eq:Yest}
\le\exp(0.3134n).
\end{equation}

For $X_0$ we directly look at all pairs $j,j'$ with $j,j'\in(n/k,n]$ with $\lcm[j,j']\le n$.
Take for example, the case $k=3$.  Then possibilities for $(a,c)$ in \eqref{eq:abc}
are $(1,1)$, $(1,2)$, and $(2,1)$.  For each $j$ we can take $j'=j$, this corresponds
to $(1,1)$.  For $j\in(n/3,n/2]$, we can also take $j'=2j$, corresponding to $(1,2)$.  
And for $j\in(2n/3,n]$ with $j$ even, we can take $j'=\frac12j$, corresponding to
$(2,1)$.  Letting $N_k(j)$ be the number of $j'$ that can correspond to $j$, we thus
have $N_k(j)=2$ for $j\in(n/3,n/2]$ and for even $j$ in $(2n/3,n]$, so that 
\[
\prod_{j\in(n/3,n]}N_3(j)\asymp 2^{n/3}.
\]
(The symbol $\asymp$ indicates the two sides are of the same magnitude up to a bounded factor.)
However, we are taking $k=30$, and this simple argument becomes more complicated, but
nevertheless can be estimated.  We have the number of assignments
for numbers $j$ in $X_0$ is
\begin{equation}
\label{eq:topend}
\le\exp(1.9115n).
\end{equation}
This estimate is arrived at as follows.   We are interested in $\prod_{j\in J_0}N_k(j)$,
which is equal to
\[
\prod_{\substack{j\in J_0\\N_k(j)>1}}N_k(j).
\]
We have $N_k(j)=1$ if and only if $j\in(n/2,n]$ and $j$ is not divisible by any prime $<k$.
Thus, up to an error of $O(1)$, the number of factors in the above product is $\nu n$, where
\[
\nu=1-1/k-(1/2)\prod_{p<k}(1-1/p).
\]
By the AM-GM inequality,
\[
\prod_{\substack{j\in J_0\\N_k(j)>1}}N_k(j)\le\left(\frac1{\nu n}\sum_{\substack{j\in J_0\\N_k(j)>1}}N_k(j)\right)^{\nu n+O(1)}.
\]
To compute the sum we refer to \eqref{eq:abc}.  By reversing the order of summation, the sum is
\[
n\sum_{\substack{a,c<k\\\gcd(a,c)=1}}\left(\frac1{ac}-\max\Big\{\frac1a,\frac1c\Big\}\right)
-\frac n2\prod_{p<k}\Big(1-\frac1p\Big)+O(1).
\]
Computing this when $k=30$ we arrive at the estimate \eqref{eq:topend}.

This leaves the contribution of numbers $j\in J_{i_0+1}$.  Note that $J_{i_0+1}$ is
a subset of $[A]$, where $A=\lceil n/(\log n)^{1/2}\rceil$.
As with the $X_i$ calculations for $i>0$, this is at most
\begin{align*}
\prod_{j\le A}(\tau(j)n/j)&\le\left(\frac1A\sum_{j\le A}\tau(j)\right)^A
\frac{n^{A}}{A!}\\
&\le(\log n)^{A}\exp\Big(A\log n-A\log A+A\Big)\\
&\le\exp\Big(\frac32 A\log\log n+A\Big).
\end{align*}
As this last estimate is of the form $e^{o(n)}$, it suffices to multiply the estimates in
\eqref{eq:Yseq}, \eqref{eq:Xest}, \eqref{eq:Yest}, and \eqref{eq:topend}, getting
that for all sufficiently large $n$, we have $\#S_\lcm(n)\le \exp(2.6071n)$.
This completes the proof.
\end{proof}

\begin{remark}
This argument gives up a fair amount in computing the contribution for $j\in X_0$,
which is the estimate \eqref{eq:topend} with $k=30$.  Another way of estimating this
count is to take some large numbers $n$ and directly compute the product
of the numbers $N_k(j)$.  It is seen that the $n$th root of this product hardly varies
as $n$ varies, and thus one can empirically arrive at a constant that is presumably
more accurate than the one in \eqref{eq:topend}.  With $k=30$, one gets in this
way the number 1.5466, which leads to the estimate $\#S_{\lcm}(n)\le\exp(2.2423n)$.
In fact, if one is prepared to reason in this way, then one can do a little better by taking
$k=100$.  This improves the numbers in \eqref{eq:Yseq}, \eqref{eq:Xest}, and \eqref{eq:Yest}
to 0.0571, 0.0807, and 0.1175, with the number in \eqref{eq:topend} moving to 1.8709,
which would give the estimate $S_\lcm(n)\le\exp(2.1262n)$.
\end{remark}

\section{Lower bounds}

Let $b$ denote a positive integer, and for $a\mid b$, let
\[
s(a,b)=\{d\mid b:d\le a\}.
\]
Further, let $p(a,b)$ denote the number of permutations $\pi$ of $s(a,b)$ such
that for each $d\in s(a,b)$, we have $\lcm[d,\pi(d)]\le a$.
Write the divisors $a$ of $b$ in increasing order: $1=a_1<a_2<\dots<a_k=b$,
where $k=\tau(b)$.  Let
\[
c(b)=\frac{\log(\tau(b)!)}{b}+\sum_{i=1}^{\tau(b)-1}\Big(\frac1{a_{i}}-\frac1{a_{i+1}}\Big)\log(p(a_{i},b)).
\]
\begin{theorem}
\label{prop:lb1}
For any positive integer $b$ we have 
\[
\#S_\lcm(n)\ge\exp((c(b)\varphi(b)/b+o(1))n)
\]
 as $n\to\infty$.
\end{theorem}
We illustrate Theorem \ref{prop:lb1} in the first interesting case: $b=2$.  Then
$p(1,2)=1$ and $p(2,2)=2$, so that $c(2)=\log(2)/2$ and the theorem asserts that
$\#S_\lcm(n)\ge\exp((\log(2)/4+o(1))n)$ as $n\to\infty$.  To see why this is true, look at sets
$\{j,2j\}$ where $j\le n/2$ and $j$ is odd.  There are $ n/4+O(1)$ of these pairs and
any permutation $\pi$ of $[n]$ for which $\pi(\{j,2j\})=\{j,2j\}$ for each $j$, and
$\pi$ otherwise acts as the identity,
is in $S_\lcm(n)$.  Since the sets $\{j,2j\}$ are pairwise disjoint, this shows that
$S_\lcm(n)$ contains at least $2^{n/4+O(1)}$ elements.  The sets are pairwise disjoint
since we are taking $j$ odd.  But a weaker condition also insures this.
Let $v_p(j)$ be the exponent on $p$ in the canonical prime factorization of $j$.
Then we take sets $\{j,2j\}$ where $j\le n/2$ and $v_2(j)$ is even.  
This insures that the sets $\{j,2j\}$ are pairwise disjoint, and now there are $n/3+O(\log n)$
pairs, leading to $\#S_\lcm(n)\ge\exp((\log(2)/3+o(1))n)$ as $n\to\infty$.

In fact, this improvement generalizes.   For a prime power $p^i$ let
\[
\alpha(p^i)=\frac{p^{i+1}-p^{i}}{p^{i+1}-1},
\]
and extend $\alpha$ as a multiplicative function on the positive integers.
Note that  $\alpha(b)$ is the density of the set of integers $j$ such that for all primes
$p\mid b$, $v_p(j)\equiv0\pmod{v_p(b)+1}$.  (Steve Fan pointed out the coincidence
to me that $\alpha(b)=b/\sigma(b)$, where $\sigma$ is the sum-of-divisors function.)
\begin{theorem}
\label{prop:lb2}
For any positive integer $b$ we have 
\[
\#S_\lcm(n)\ge\exp((c(b)\alpha(b)+o(1))n)
\]
 as $n\to\infty$.
\end{theorem}
\begin{proof}
For $1\le i\le\tau(b)-1$ consider the intervals $I_i:=(n/a_{i+1},n/a_{i}]$ and $I_{\tau(b)}=(0,n/b]$.
For $j\in I_i$ with $v_p(j)\equiv0\pmod{v_p(b)+1}$ for each prime $p\mid b$, we have
the set $T(i,j):=\{dj:d\in s(a_i,b)\}$ as a subset of $[n]$.  Moreover, the sets $T(i,j)$ are pairwise 
disjoint for all pairs $i,j$ with $j\in I_i$ and $v_p(j)\equiv0\pmod{v_p(b)+1}$ for
each prime $p\mid b$.  For any permutation $\pi$ of $[n]$, if $\pi(T(i,j))=T(i,j)$
for all $i,j$ and otherwise $\pi$ acts as the identity, we have $\pi\in S_\lcm(n)$.

For a given value of $i<\tau(b)$ there are
$\sim(1/a_i-1/a_{i+1})\alpha(b)n$ values of $j$, and for $i=\tau(b)$, there are 
$\sim n\alpha(b)/a_{\tau(b)}$ values of $j$.  We conclude that $\#S_\lcm(n)$ is
at least 
\[
(p(a_\tau(b),b))^{(1+o(1))n\alpha(b)/a_{\tau(b)}}\prod_{i=1}^{\tau(b)-1}p(a_i,b)^{(1+o(1))n\alpha(b)(1/a_i-1/a_{i-1})}.
\]
Since $a_{\tau(b)}=b$ and $p(b,b)=\tau(b)!$, the result follows.
\end{proof}

We have an analogous result for $S_\ddiv(n)$.
Let $p_\dd(a,b)$ denote the number of permutations $\pi$ of $s(a,b)$ such
that for each $d\in s(a,b)$, we have $d\mid\pi(d)$ or $\pi(d)\mid d$.
Let
\[
c_\dd(b)=\frac{\log(p_\dd(b,b))}{b}+\sum_{i=1}^{\tau(b)-1}\Big(\frac1{a_{i}}-\frac1{a_{i+1}}\Big)\log(p_\dd(a_{i},b)).
\]
\begin{corollary}
\label{cor:lb2}
For any positive integer $b$ we have 
\[
\#S_\ddiv(n)\ge\exp((c_\dd(b)\alpha(b)+o(1))n)
\]
 as $n\to\infty$.
\end{corollary}

\begin{footnotesize}
\begin{table}[]
\caption{Some values of $c(b)\alpha(b)$ and $c_\dd(b)\alpha(b)$ to 6 places with their
exponentials rounded down to 4 places.}
\label{Ta:A}
\begin{tabular}{|rrrrrr|} \hline
$b$&$c(b)\alpha(b)$ & $e^{c(b)\alpha(b)}$&$c_\dd(b)\alpha(b)$&$e^{c_\dd(b)\alpha(b)}$&\\ \hline
4&.354987&1.4261&.354987&1.4261&\\
12&.536243&1.7095&.479872&1.6158&\\
24&.602065&1.8258&.542689&1.7206&\\
48&.638300&1.8932&.578122&1.7826&\\
60&.646856&1.9095&.552061&1.7368&\\
96&.658201&1.9313&.597849&1.8182&\\
120&.707611&2.0291&.610358&1.8410&\\
144&.704928&2.0237&.631752&1.8809&\\
210&.600981&1.8239&.496559&1.6430&\\
240&.740127&2.0962&.648821&1.9132&\\
288&.723607&2.0618&.650371&1.9162&\\
420&.716176&2.0465&.597383&1.8173&\\
480&.757765&2.1335&.660864&1.9364&\\
\hline
\end{tabular}
\end{table}
\end{footnotesize}
 So, by Table \ref{Ta:A} and Theorem \ref{prop:lb2} with $b=480$, we have
$P_\lcm(n)\ge 2.1335^n$ for all large values of $n$ and by Corollary \ref{cor:lb2}
we have $P_\ddiv(n)\ge 1.9364^n$ for all large $n$.

\section{Comparing $\#S_\ddiv(n)$ and $\#S_\lcm(n)$}
Let
\[
R(n)=\#S_\lcm(n)/\#S_\ddiv(n).
\]
It appears from a glance at Table \ref{Ta:counts} that $R(n)$ grows at least geometrically.
Here we prove this.
\begin{theorem}
\label{thm:ratio}
There is a constant $c>1$ such that for all large values of $n$ we have $R(n)>c^n$.
\end{theorem}
\begin{proof}
Let $A$ denote the set of integers $a$ with $n/7<a\le n/6$ and with $a$ not divisible by any prime $<10^4$.   Note that
\[
\#A\ge \frac n{42}\prod_{p<10^4}\Big(1-\frac1p\Big)+O(1),
\]
so that $\#A>14n/10^4$ for all large $n$.
 For $a\in A$, let 
 \[
 B_a=\{a,2a,3a,4a,5a,6a\}.
 \]
   Any divisor of a member of $B_a$ that is
 not itself in $B_a$ must be $<(n/6)/10^4$.  Since clearly each member of $B_a$ has no
 multiple in $[n]$ that is not in $B_a$, we have that each $\pi\in S_\ddiv(n)$ has
 $\pi(B_a)\ne B_a$ for at most $n/10^4$ values of $a\in A$.  We conclude that
 each $\pi\in S_\ddiv(n)$ has $\pi(B_a)=B_a$ for at least $13n/10^4$ values of $a\in A$.
 
 Let $\pi\in S_\ddiv(n)$ with $\pi(B_a)=B_a$, and let $\pi_0$ be $\pi$ restricted to
 $[n]\setminus B_a$.  There are exactly 36 permutations $\pi\in S_\ddiv(n)$ which
 give rise to the same $\pi_0$ corresponding to the 36 permutations $\sigma$ of $B_a$
 with each $ja\mid\sigma(ja)$ or $\sigma(ja)\mid ja$ (since $\#S_\ddiv(6)=36$).  
 However, for a given $\pi_0$ here, there are exactly 56 permutations $\pi\in S_\lcm(n)$
 with $\pi$ restricted to $[n]\setminus B_a$ equal to $\pi_0$ (since $\#S_\lcm(6)=56$). 
 
It thus follows that
 \[
 R(n)>(56/36)^{13n/10^4},
 \]
 so the theorem is proved with $c=(56/36)^{13/10^4}>1.00057$.
 \end{proof}
 
 \section{An upper bound for anti-coprime permutations}
 
 Let $A(n)$ denote the number of anti-coprime permutations $\pi$ of $[n]$.
 We prove the following theorem.
 \begin{theorem}
 \label{thm:anti}
 We have $A(n)=n!/(\log n)^{(e^{-\gamma}+o(1))n}$ as $n\to\infty$.
 \end{theorem}
 \begin{proof}
 In light of the lower bound from \cite{P22}, it suffices to prove that
 $A(n)\le n!/(\log n)^{(e^{-\gamma}+o(1))n}$ as $n\to\infty$.

Let $\psi(n)\to\infty$ arbitrarily slowly, but with $\psi(n)=o(\log\log n)$.  
Let $\alpha=1+1/(\psi(n))^{1/2}$, so that $\alpha\to1^+$ as $n\to\infty$.
 Let the integer variable $i$ satisfy
 \begin{equation}
 \label{eq:Ii}
  \psi(n)<i<\log\log n/\log\alpha-\psi(n).
 \end{equation}
 Thus, $e^{\alpha^i}\to\infty$ and $e^{\alpha^i}=n^{o(1)}$ as $n\to\infty$.
 For each $i$ satisfying \eqref{eq:Ii}, let
 \[
I_i:=(e^{\alpha^{i-1}},e^{\alpha^i}],\quad J_i:=\{j\in[n]:P^-(j)\in I_i\},
 \]
 where $P^-(j)$ is the least prime factor of $j$.
 By the Fundamental Lemma of the Sieve and Mertens' theorem, we have
 \begin{equation}
 \label{eq:Ji}
 \#J_i\sim \frac n{e^\gamma \alpha^{i-1}}-\frac n{e^\gamma\alpha^i}=\frac{n(\alpha-1)}{e^\gamma\alpha^i}
 \end{equation}
 uniformly in $i$ satisfying \eqref{eq:Ii}, as $n\to\infty$.
 
 For $j\in[n]$ the number of $j'\in[n]$ with $\gcd(j',j)>1$ is $\le\sum_{p\mid j}n/p$.
 (This is a poor bound for most integers $j$, but fairly accurate for most $j$'s without
 small prime factors, as is the case for members of $J_i$.)
 Thus, the total number of assignments for the numbers $j\in J_i$ in an anti-coprime
 permutation of $[n]$ is
 \[
 \le\prod_{j\in J_i}\left(n\sum_{p\mid j}\frac1p\right)\le
 \left(\frac n{\#J_i}\sum_{j\in J_i}\sum_{p\mid j}\frac1p\right)^{\#J_i},
 \]
 by the AM-GM inequality.  The double sum here is
 \[
 \sum_{p>e^{\alpha^{i-1}}}\sum_{\substack{j\in J_i\\p\mid j}}\frac1p
 \ll n\sum_{p>e^{\alpha^{i-1}}}\frac1{p^2\alpha^i}\ll\frac n{e^{\alpha^{i-1}}\alpha^{2i}},
 \]
 uniformly, using an upper bound for the sieve.  Thus, for $n$ large,
 \[
 \frac n{\#J_i}\sum_{j\in J_i}\sum_{p\mid j}\frac1p\le\frac n{(\alpha-1)e^{\alpha^{i-1}}}
 \le\frac n{e^{\alpha^{i-2}}}.
 \]
 Hence, using \eqref{eq:Ji},
 \begin{align*}
 \prod_{j\in J_i}\left(n\sum_{p\mid j}\frac1p\right)&\le
 \left(\frac n{e^{\alpha^{i-1}}}\right)^{\#J_i}
  =\exp(\#J_i\log n-\#J_i\alpha^{i-2})\\
  &=\exp(\#J_i\log n-{(1+o(1))n(\alpha-1)/e^\gamma})
  \end{align*}
  uniformly.

Let $N$ be the number of $j\in[n]$ not in any $J_i$, so that $N=n-\sum_i\#J_i$.
Further the number of values of $i$ is at most $\log\log n/\log\alpha-2\psi(n)$.
After assignments have been made for the values of $j \in \cup J_i$, there
are at most $N!$ remaining assignments for $j\notin\cup J_i$ and multiplying
this by the product of the previous estimate for all $i$ is at most
\[
\exp\left(n\log n-(1+o(1))\frac{n(\alpha-1)}{e^\gamma}\Big(\frac{\log\log n}{\log\alpha}-2\psi(n)\Big)\right).
\]
 It remains to note that $(\alpha-1)/\log\alpha\sim1$ as $n\to\infty$, which completes
 the proof of the theorem.
 \end{proof}
 
 \section*{Dedication and Acknowledgments}
 Eduard Wirsing is not primarily known for his work in combinatorial number theory, yet one
 of his papers that influenced me a great deal is his joint work with Hornfeck, later improved
 on his own (see \cite{W}), on the distribution of integers $n$ with $\sigma(n)/n=\alpha$, for a fixed 
 rational number $\alpha$, where $\sigma$ is the sum-of-divisors function.  In a survey I wrote
 with S\'ark\"ozy \cite{PS} on combinatorial number theory, we singled out this particular work for 
 being a quintessential exemplar of the genre.  It is in this spirit that I offer this note on combinatorial
 number theory in remembrance of Eduard Wirsing.
 
  I wish to thank \v{S}ar\={u}nas Burdulis, Steve Fan,
  Mitsuo Kobayashi, Jared Lichtman, and Andreas Weingartner
 for help in various ways.


\begin{thebibliography}{99}

\bibitem{BB}
E. Berkove and M. Brilleslyper, Subgraphs of coprime graphs on sets of consecutive
integers, Integers {\bf 22} (2022), \#A47.

\bibitem{EFH}
P. Erd\H os, R. Freud, and N. Hegyv\'ari,
Arithmetical properties of permutations of integers,
Acta Math. Hungar. {\bf 41} (1983), 169--176.

\bibitem{M}
N. McNew,
Counting primitive subsets and other statistics of the divisor graph of $\{1,2,\dots,n\}$,
 European J. Combin. {\bf92} (2021), Paper No.\ 103237, 20 pp.

\bibitem{oeis}
Online Encyclopedia of Integer Sequences, A320843.

\bibitem{P2}
C. Pomerance, On the longest simple path in the divisor graph,
Proceedings of the fourteenth Southeastern conference on combinatorics, graph theory and computing (Boca Raton, Fla., 1983).
Congr. Numer. {\bf40} (1983), 291--304.


\bibitem{P22}
C. Pomerance, Coprime permutations, arXiv:2203.03085 [math.NT], 20 pp.

\bibitem{PS}
C. Pomerance and A. S\'ark\"ozy, Combinatorial number theory,
in Handbook of Combinatorics, R. L. Graham, M. Gr\"otschel, L. Lov\'asz, eds., Elsevier Science B.V., 1995, pp. 967--1018.

\bibitem{SS}
A. Sah and M. Sawhney, Enumerating coprime permutations, arXiv:2203.06268 [math.NT], 11 pp.

\bibitem{S}
E. Saias, Etude du graphe divisoriel 5, arXiv:2107.03855 [math.NT], 39 pp.

\bibitem{W}
E. Wirsing, Bemerkung zu der Arbeit \"uber volkommene Zahlen,
Math. Ann. {\bf137} (1959), 316--318.

\end{thebibliography}
\end{document}